\documentclass{article}

\usepackage[utf8]{inputenc}
\usepackage[russian]{babel}
\usepackage{graphicx,amsfonts}

\begin{document}
\def\newstr{\par\noindent}
\def\ms{\medskip\newstr}
\def\stat#1#2{{\bf #1}: {\it #2}}
\def\ra{\longrightarrow}
\def\ovl{\overline}
\def\op{\oplus} \def\om{\ominus}

\renewcommand{\labelenumi}{\theenumi)}
\def\f{\varphi}
\def\A{\mathfrak A}

\ \\ \ \\
\centerline{\huge\bf Another Solution to the Thue Problem}

\smallskip\smallskip
\centerline{\huge\bf of Non-Repeating Words}

\medskip\ms\centerline{\bf\large Boris Zolotov (Saint-Petersburg, Russia)}


\ \\ [3.5cm]
\newstr\centerline{\Large\bf Abstract}

\begin{quote}
\newstr In this work we consider morphisms that preserve well-known non-repeating properties: squarefreeness, cubefreeness, overlap-freeness and weak squarefreeness. Up to the present moment only the morphisms preserving three out of four non-repeating properties have been known. The problem of the existence of weakly squarefree morphisms was open.

\ms The essential result of this work is the positive solution to this problem. An example of the morphism preserving all four properties is provided. Also, it is proved that there are no morphisms with the same properties and a lower rank.
\end{quote}

\vfill\eject


\ \ 
\section{Introduction. The Problem}

The problem was raised at the beginning of XX century by a Norwegian mathematician Axel Thue. He wondered if it was possible to construct an infinite word over the least possible alphabet such as no finite factors of it repeat twice or three times. Such words were called squarefree or cubefree, respectively. Thue [1] proved that such words exist and provided corresponding examples.

\ms An English mathematician John Leech [3] worked on the same problem. He found an example of an infinite squarefree word, and this word was constructed using an uniform morphism with the rank 13. An uniform morphism is a morphism such as all the letters have images of the same length, and a rank of the morphism is the length of the images of the letters.

\ms In 1982 Max Crochemore [4] found an algorithmically verifiable criterion of squarefreeness of a morphism.

\ms In this work new results are obtained that allow us to optimize Thue's and Leech's results and construct new sets of infinite non-repeating words.


\medskip\medskip
\section{Known Results}
Let us provide Thue's, Leech's and Crochemore's Results.

\begin{enumerate}
  \item The first example of an infinite cubefree word (Thue [1], 1906).
    \newstr Let $\A=\{0;1\}$; $A_0=1$; $\f: 1 \ra 10, 0 \ra 01$;
    \subitem $A_0 = 1$
    \subitem $\f(A_0) = 1\ 0$
    \subitem $\f^2(A_0) = 10\ 01$
    \subitem $\f^3(A_0) = 1001\ 0110$
    \subitem $.\ .\ .\ .$
    \subitem $\f ^\infty (A_0) = 1001011001101001...$ is an infinite cubefree word, which is also called the Thue-Morse sequence.

\ms For the morphism $\f$ some other properties were also proved. Thue [1] showed that the fixed point $\f ^\infty (A_0)$ did not contain any factors of the form $aXaXa$, which are called {\it overlaps}. Furthermore, $\f(W\!)$ is cubefree when $W\!$ is cubefree, and morphisms with this property are called {\it cubefree}. Then, the image of any overlap-free word is overlap-free, so $\f$ is also {\it overlap-free}. Morphisms with such properties are called {\it Thue morphisms}.

  \item An example of morphism that generates an infinite squarefree word. (Thue [2], 1912).
    \newstr Let $\A=\{1;2;3\}$; $A_0=1$; $\f: 1 \ra 12312, 2 \ra 131232, 3 \ra 1323132$.

\vfill\eject

  \item An uniform morphism that generates an infinite squarefree word (Leech [3], 1957).
   \newstr Let $\A=\{1;2;3\}$; $A_0=1$;\\
    $\f: 1 \ra 1232132312321, 2 \ra 2313213123132, 3 \ra 3121321231213$.

  \item The characterization of squarefreeness of a morphism (Crochemore [4], 1982).
    \newstr A morphism defined on a three letter alphabet is squarefree iff the images of the squarefree words of length 5 are squarefree.

  \item The characterization of squarefreeness of an uniform morphism (Crochemore [4], 1982).
    \newstr An uniform morphism defined on a three letter alphabet is squarefree iff the images of the squarefree words of length 3 are squarefree.

\end{enumerate}


\medskip\medskip
\section{Definitions}

{\bf Alphabet, Word, Morphism}

\smallskip\newstr An {\it alphabet} is a finite set, and its elements are called {\it letters}. A {\it word upon the alphabet}  is a sequence of this alphabet's letters. The number of letters in the word $W\!$ is called the {\it length} of $W\!$ and denoted by $|W\!|$. The $i$--th letter in the word $W\!$ is denoted by $W\![i]$.

\ms The set of all the words upon the alphabet $\A$ is denoted by $\ovl{\A}$. The concatenation operation makes it be a free monoid with an empty word as a neutral element.

\ms A {\it morphism} is a function $\f \colon \ovl{\A} \mapsto \ovl{\A}$ that preserves the concatenation: \smallskip\\
\centerline{$\forall\,u,v \in \ovl{\A} \quad \f(uv) = \f(u)\f(v).$}

\ms Let $L$ be a positive integer. A morphism $\f$ is called {\it $L$--uniform} if for every $a \in \A$ we have $|\f(a)|=L$. In this case, $L$ is called the {\it rank} of the morphism $\f$.

\medskip\ms{\bf Morphism and fixed point}

\smallskip\newstr As long as $\f$ is an automorphism of a free monoid $\ovl{\A}$, it is sufficient to define images of the elements of $\A$ to entirely describe $\f$.
\newstr An example:
\newstr\qquad Let $\A=\{1;2;3\}$ and let $X=1231$ be a word upon $\A$.
\newstr\qquad $\f: 1 \ra 1231, 2 \ra 3112, 3 \ra 1132$.
\newstr\qquad Then $\f(1231) = \f(1)\f(2)\f(3)\f(1) = 1231\,3112\,1132\,1231$.

\ms When $r$--uniform morphism $\f$ is applied to a word $W\!$, the length of $\f(W\!)$ is divisible by $r$. So, we can define {\it canonical fragments}, which are non-crossing factors of length $r$ in $\f(W\!)$ that cover $\f(W\!)$ and are images of the letters of $W\!$.

\ms Let us consider an iterational sequence of words $A_n$ such as $A_{n+1} = \f(A_n)$. If it is known that for any integer $n$ $A_{n+1}=A_nV_n$ then there exists an infinite word $A$ such as:
\begin{enumerate}
  \item for every integer $n$ it is true that $A=A_nW\!_n$;
  \item $\f(A)=A$.
\end{enumerate}
\newstr $A$ is called the {\it fixed point} of the morphism $\f$ and denoted by $\f^\infty (A_1)$.


\vfill\eject
\newstr{\bf Increasing and Decreasing, Cyclic Morphisms}

\smallskip\newstr Let us consider the ternary alphabet $\A=\{1;2;3\}$. We will define {\it operations} `$\op 1$' and `$\om 1$' upon it, such as: \\
\centerline{ $\begin{array}{l}
    1 \op 1=2 \\
    2 \op 1=3 \\
    3 \op 1=1 \\
\end{array}$ \qquad
$\begin{array}{l}
    1 \om 1=3 \\
    2 \om 1=1 \\
    3 \om 1=2 \\
\end{array}$.}

\ms Given a word $T$ upon a ternary alphabet, we can define words $T \op 1$ and $T \om 1$.

\smallskip\newstr $T \op 1$ is a word such that \\
\centerline{$\forall i = \ovl{1,|T|}$ \quad $(T \op 1)[i] = (T[i]) \op 1$;}
\newstr $T \om 1$ is a word such that \\
\centerline{$\forall i = \ovl{1,|T|}$ \quad $(T \om 1)[i] = (T[i]) \om 1$.}

\ms $r$--uniform morphism $\f$ is called {\it cyclic} if \\
\centerline{$\forall\, p \in \A$ \quad $\f(p \op 1) = (\f(p)) \op 1$.}

\ms A word $W\!$ of length $l$ is called {\it increasing} if \\
\centerline{$\forall\, i = \ovl{1,l-1}$ \quad $W\![i+1]=W\![i] \op 1$.}

\ms A word $W\!$ of length $l$ is called {\it decreasing} if \\
\centerline{$\forall\, i = \ovl{1,l-1}$ \quad $W\![i+1]=W\![i] \om 1$.}

\medskip\ms{\bf Properties of words and morphisms}

\smallskip\newstr A {\bfseries\itshape square} is a word of the form $X\!X$, where $X$ is not empty. Example --- $123\,123$.
\smallskip\newstr A word is called {\it squarefree} if it does not contain any squares.
\smallskip\newstr A morphism is called {\it squarefree} if image of every squarefree word is squarefree. A morphism is called {\it $k$--squarefree} if image of every squarefree word of length $k$ is squarefree.

\smallskip\ms A {\bfseries\itshape cube} is a word of the form $X\!X\!X$, where $X$ is not empty. Example --- $123\,123\,123$.
\smallskip\newstr A word is called {\it cubefree} if it does not contain any cubes.
\smallskip\newstr A morphism is called {\it cubefree} if the image of every cubefree word is cubefree.

\smallskip\ms An {\bfseries\itshape overlap} is a word of the form $aX\!aX\!a$, $X$ can be empty. Example --- $2\,123\,2\,123\,2$.
\smallskip\newstr A word is called {\it overlap-free} if it does not contain any overlaps.
\smallskip\newstr A morphism is called {\it overlap-free} if the image of every overlap-free word is overlap-free.

\smallskip\ms A {\bfseries\itshape weak square} is a word of the form $aX\!X\!a$, $X$ can be empty. Example --- $3\,123\,123\,3$.
\smallskip\newstr A word is called {\it weakly squarefree} if it does not contain any weak squares.
\smallskip\newstr A morphism is called {\it weakly squarefree} if the image of every weakly squarefree word is weakly squarefree.

\smallskip\ms A morphism is called {\it a Thue morphism} if it is cubefree, overlap-free and generates a fixed point $\f ^\infty (a)$, where $a$ is a letter.

\vfill\eject


\section{Objectives}
In this research we are addressing the following questions:

\begin{enumerate}
  \item Is Leech's result optimal? Are there any uniform squarefree morphisms with a lower rank?
  \item Are there any uniform Thue morphisms over the ternary alphabet?
  \item Are there any weakly squarefree Thue morphisms over the ternary alphabet?
  \item How is weak squarefreeness connected with other properties, for example, with usual squarefreeness?
  \item Which combinations of letters are necessary in squarefree words? What is the necessary length for these combinations?
\end{enumerate}


\section{Main Results}
The main results of this work are the following:

\ms{\bf Newly Created Morphisms}
\newstr The following morphisms were considered, and their properties were studied:
\begin{itemize}
  \item[\bf 1)] A morphism  $\f: 1 \ra 12321, 2 \ra 23132, 3 \ra 31213$ is a cubefree and weakly squarefree morphism over the ternary alphabet. It also generates a fixed point.
  \item[\bf 2)] A morphism  $\f: 1 \ra 1221, 2 \ra 2332, 3 \ra 3113$ is a Thue morphism over the ternary alphabet.
  \item[\bf 3)] A morphism  $\f: 1 \ra 121, 2 \ra 232, 3 \ra 313$ is a cubefree morphism over the ternary alphabet. It also generates a fixed point.
\end{itemize}

\ms{\bf Low-Rank Morphisms}
\newstr The following theorems were proved for the low-rank morphisms:

\begin{itemize}
  \item[\bf 4)] Over the ternary alphabet there are no uniform weakly squarefree Thue morphisms with a rank lower than 5.
  \item[\bf 5)] Over the binary alphabet there are no weakly squarefree Thue morphisms.
\end{itemize}

\ms{\bf Weak Squarefreeness}
A theorem on connection between non-repeating properties was proved:

\begin{itemize}
  \item[\bf 6)] An uniform and cyclic squarefree morphism is weakly-squarefree.
\end{itemize}

\ms{\bf Leech's Morphism and Optimal Ranks}
\newstr The following results were obtained:
\begin{itemize}
  \item[\bf 7)] There are exactly 144 uniform squarefree morphisms of rank 11; there are no uniform squarefree morphisms with a lower rank.
  \item[\bf 8)] The Leech's morphism is cyclic; there are no cyclic squarefree morphisms with a rank lower than 13.
  \item[\bf 9)] The Leech's morphism is squarefree, cubefree, overlap-free and weakly squarefree. There are no morphisms with such properties and a lower rank.
\vfill\eject
\end{itemize}

\ms{\bf Letter combinations in squarefree words}
\newstr Letter combinations over the ternary were considered, and the following results were obtained:

\begin{itemize}
  \item[\bf 10)] Each squarefree word of length more than 13 over the ternary alphabet contains all the two-letter combinations of distinct letters.
  \item[\bf 11)] Each squarefree word of length more than 36 over the ternary alphabet contains all the three-letter combinations of distinct letters (of the form $abc$).
\end{itemize}

\ms{\bf Summary}

\smallskip\newstr All the questions addressed were answered. What was found is the lower bound of the rank of uniform squarefree morphisms over the ternary alphabet. All the morphisms having such a rank were provided.

\ms Furthermore, it was shown that there exist weakly squarefree Thue morphisms over the ternary alphabet. The lower bound of the rank of such morphisms was provided. Also, we provided low-rank morphisms over the ternary alphabet that have only a few of the considered properties. It proves that the properties are independent from each other.

\ms On the other hand, it was proved that cyclicness and squarefreeness result in weak squarefreeness. Also, the properties of the letter combinations were examined.

\ms There remains the open problem of the existence of uniform weakly squarefree Thue morphisms of rank between 5 and 13.


\medskip\medskip\medskip
\section{Proofs of the main results. \\ Newly created morphisms}

\stat{Lemma (0)}{A morphism $\f$ has a fixed point $\f ^\infty (u)$ (where $u$ is a letter) iff there exists a word $V$ such as $\f(u) = uV$.}

\ms We will use the induction method to prove the lemma. Let us construct a sequence of words $A_n$ such as $A_0=u$ and $A_{k+1} = \f(a_k)$. It is necessary to show that for any $k \in \mathbb N_0$ there exists a word $W\!$ such as $A_{k+1} = A_kW\!$.

\ms For $k=0$ the statement is true, because it is given in the wording. Now let us assume that for an integer $k$ we have $A_{k+1} = A_kW\!$.

\ms Then $A_{k+2} = \f(A_{k+1}) = \f(A_kW\!) = \f(A_k)\f(W\!) = A_{k+1}W'\!$. As a result, we obtain that each word in the sequence $A_n$ is a prefix of the next one. Consequently, the morphism $\f$ has a fixed point $\f ^\infty (u)$.

\medskip\ms\stat{Theorem 1}{The morphism $\f: 1 \ra 12321, 2 \ra 23132, 3 \ra 31213$ is a cubefree morphism over the ternary alphabet.}

\begin{center}
\includegraphics[width=4.7cm]{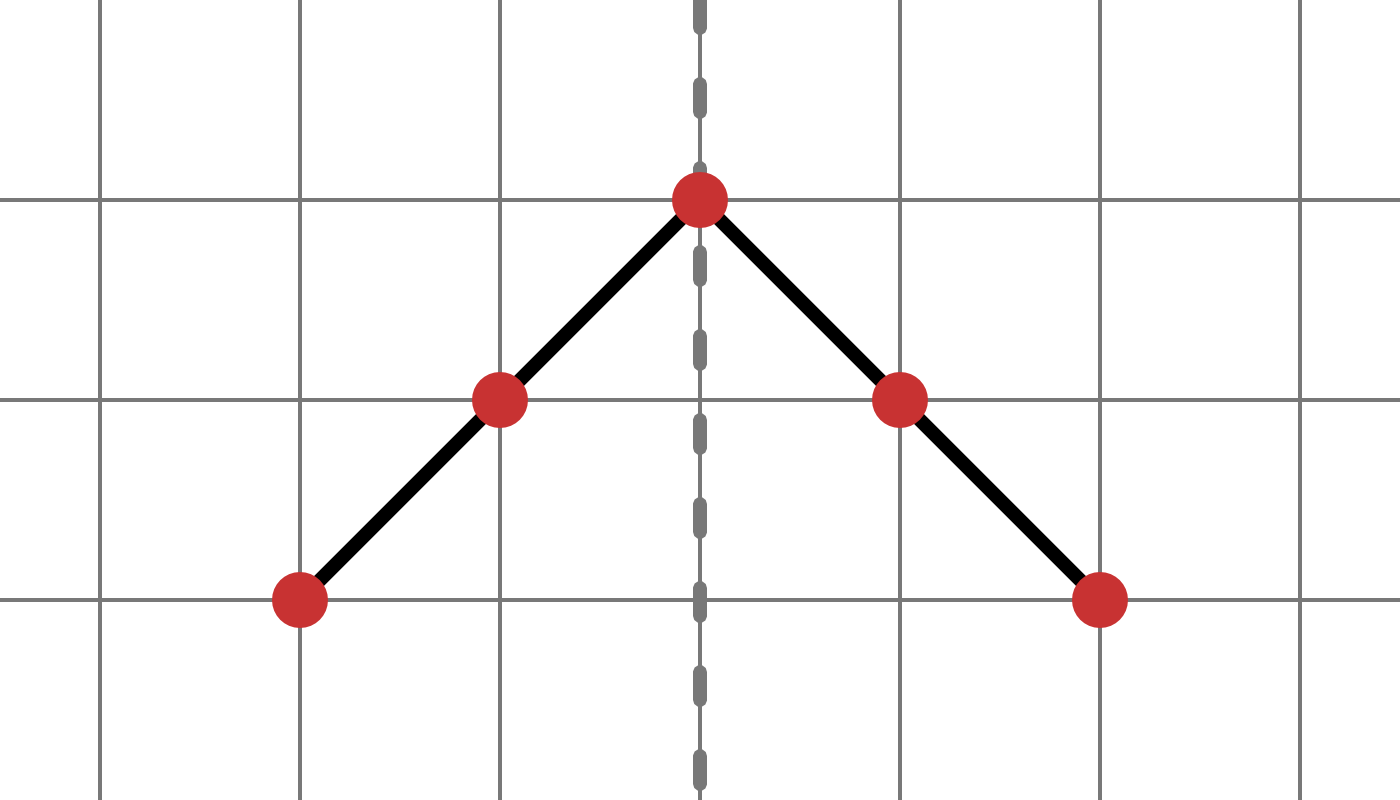}
{\\ \small\it a canonical fragment of the morphism $\f$, \\
increasing and decreasing factors of it}
\end{center}

\vfill\eject
\ms We need a series of lemmas to prove this theorem.

\medskip\ms\stat{Lemma (1)}{When morphism $\f$ is applied to a cubefree word, the image of this word does not contain any cubes of letters (of the form $aaa$).}

\ms Let the image of a cubefree word have a cube of the letter $aaa$. Then two of these three letters belong to one canonical fragment. On the other hand, it is easy to see that inside the canonical fragments the letters that are next to each other are different.

\medskip\ms\stat{Lemma (2)}{When morphism $\f$ is applied to a cubefree word, the image of this word does not contain any cubes of the form $X\!X\!X$, where $|X|=2$.}

\ms Let us notice that in the canonical fragment $a_1a_2a_3a_4a_5$ the following equations are true:
\begin{center}
$a_1=a_5$; \\
$a_2=a_4$; \\
$a_2=a_1 \op 1$; \\
$a_3=a_2 \op 1$.
\end{center}

\ms Now let us denote the cube $X\!X\!X$ by $c_1c_2c_3c_4c_5c_6$. Two canonical fragments are sufficient to cover the cube; we will denote these canonical fragments by $a_1a_2a_3a_4a_5$ and $b_1b_2b_3b_4b_5$.

\ms Finally, we are to consider all the possible cases of the position of the letter $c_1$ inside the first canonical fragment.

\begin{enumerate}
\item $c_1=a_1$. Then  $c_3=(c_1 \op 1) \op 1=c_1 \om 1$. This contradicts the fact that $c_3=c_1$.
\item $c_1=a_2$. Then the word $c_1c_2$ is increasing, whereas $c_3c_4$ is decreasing. This contradicts the fact that these words are equal.
\item $c_1=a_3$. Then the word $c_1c_2$ is decreasing, whereas $c_5c_6$ is increasing. This contradicts the fact that these words are equal.
\item $c_1=a_4$. Then the word $c_1c_2$ is decreasing, whereas $c_3c_4$ is increasing. This contradicts the fact that these words are equal.
\item $c_1=a_5$. Then the word $c_3c_4$ is increasing, whereas $c_5c_6$ is decreasing. This contradicts the fact that these words are equal.
\end{enumerate}

\ms As a result, the cube of length 6 can not begin with any letter inside the canonical fragment. The Lemma is proved.

\medskip\ms\stat{Lemma (3)}{When morphism $\f$ is applied to a cubefree word, the image of this word does not contain any cubes of the form $X\!X\!X$, where $|X|=3$.}

\ms Several cases are possible of how the cube $X\!X\!X$ lies inside the sequence of canonical fragments. The first one is when the second occurence of $X$ consists of the 2-nd, the 3-rd and the 4-th letter of a canonical fragment. Then the word made of two last letters of $X$ is increasing and decreasing at the same time; this case is impossible.

\ms In other cases the following fact is easy to check: two words $X$ out of three must lie inside canonical fragments without crossing the borders. These two words will begin at different positions inside the canonical fragments, because the length of the cube is not divisible by 5.

\ms Now let us notice that all the distinct factors of length 3 inside the canonical fragment are different if increasing\,/\,decreasing is considered. This fact makes the existence of the cube $X\!X\!X$ impossible, if $|X|=3$.
\vfill\eject


\medskip\ms\stat{Lemma (4)}{When morphism $\f$ is applied to a cubefree word, the image of this word does not contain any cubes of the form $X\!X\!X$, where $|X|=4$.}

\ms To prove this Lemma it is also needed to consider five cases of the position of the first letter of the cube inside the canonical fragment. In each case it is easy to obtain that a factor of $X$ is increasing and decreasing at the same time.

\medskip\ms\stat{Lemma (5)}{When morphism $\f$ is applied to a cubefree word, the image of this word does not contain any cubes of the form $X\!X\!X$, where $|X|$ is divisible by 5.}

\ms Let us notice that the morpism $\f$ is cyclic. For a cyclic morphism it is true that the canonical fragment can be determined by a letter inside it and the position of this letter inside the fragment.

\ms Furthermore, the corresponding letters in the words $X$ are located on the same positions inside canonical fragments. It leads to the fact that the corresponding canonical fragments are identical. Finally, we can say that the original word contained a cube, which contradicts the wording.

\medskip\ms\stat{Lemma (6)}{A sequence of canonical fragments of the morphism $\f$ does not contain any squares $X\!X$ such that $|X|$ is greater than 5 and not divisible by 5.}

\ms Let such a square exist inside the sequence of fragments. Then each word $X$ contains a border between two canonical fragments next to each other. Furthermore, if we compare these two words $X$, we will see that the borders of the canonical fragments are shifted in the second one.

\ms On the other hand, it can be noticed that any shift of the canonical fragment of the morphism $\f$ has increases\,/\,decreases differently from the original one. Finally, it means that one factor of $X$ has to be increasing and decreasing at the same time, which is impossible.

\ms The Lemma is now proved.

\medskip\ms Now let the image of a cubefree word contain a cube. Lemmas 1--6 consider all the cases of the length of $X$, and for each case it is proved that cubes with this length do not exist.

\ms The theorem is proved.


\medskip\ms\stat{Theorem 2}{The morphism $\f: 1 \ra 12321, 2 \ra 23132, 3 \ra 31213$ is a weakly squarefree morphism over the ternary alphabet.}

\medskip\ms\stat{Lemma (7)}{When morphism $\f$ is applied to a weakly squarefree word, the image of this word does not contain any squares of a letter.}

\ms Let a square of a letter exist inside the image of a weakly squarefree word. Two identical successive letters can not lie inside one canonical fragment, because inside the canonical fragment all successive letters are different. It means that equal letters lie on the edges of two successive canonical fragments. As a result, these fragments are the same, and in the original word there was a square of a letter, which is also a weak square. This contradicts the wording.

\medskip\ms\stat{Lemma (8)}{When morphism $\f$ is applied to a weakly squarefree word, the image of this word does not contain any weak squares of the form $abba$, where $a$, $b$ are letters.}

\ms Let us notice that the factor $bb$ can not lie inside a canonical fragment, because successive letters inside the canonical fragment are different. It means that the weak square $abba$ lies on the border of two canonical fragments, one of which ends with $b$, and the other one starts with $b$.

\vfill\eject
\ms Now let us note that the only canonical fragment can start with $b$, and it is $\f(b)$. Also, $\f(b)$ is the only canonical fragment that ends with $b$. It means that the original word contained the weak square $bb$, which contradicts the wording.

\medskip\ms\stat{Lemma (9)}{When morphism $\f$ is applied to a weakly squarefree word, the image of this word does not contain any weak squares of the form $aX\!Xa$, where $|X|=2$.}

\ms We will prove that image of the weakly squarefree word can not contain a square $X\!X$, where $|X|=2$. If such a square lies inside the image, then two canonical fragments are sufficient to cover it. Let us denote the first canonical fragment by $a_1a_3a_3a_4a_5$ and consider the cases of the position of the first letter of the square inside it.

\begin{enumerate}
  \item If the first letter of the square is the first or the second letter of the canonical fragment, then the word $X$ is increasing and decreasing at the same time, which is impossible.
  \item If the first letter of the square is the third letter of the canonical fragment, then the first letters of two words $X$ are different, which is impossible.
  \item If the first letter of the square is the fourth letter of the canonical fragment, then the word $X$ is increasing and decreasing at the same time, which is impossible.
  \item If the first letter of the square is the fifth letter of the canonical fragment, then the last letters of two words $X$ are different, which is impossible.
\end{enumerate}

\ms Now the Lemma is proved.

\medskip\ms\stat{Lemma (10)}{When morphism $\f$ is applied to a weakly squarefree word, the image of this word does not contain any weak squares of the form $aX\!Xa$, where $|X|=3$.}

\ms We will prove that image of the weakly squarefree word can not contain a square $X\!X$, where $|X|=3$. If such a square exists then the first or the last letters of two words $X$ lie in one canonical fragment.

\ms Now let us notice that two letters inside one canonical fragment are not identical if the difference between their positions is equal to 3. It means that the first or the last letters of $X$ are different, which is impossible.

\medskip\ms\stat{Lemma (11)}{When morphism $\f$ is applied to a weakly squarefree word, the image of this word does not contain any weak squares of the form $aX\!Xa$, where $|X|=4$.}

\ms In the beginning, let us notice that the image of a weakly squarefree word can contain a square $X\!X$, $|X|=4$. For example, $\f(212) = 231\,3212\,3212\,3132$. This square is not weak, though.

\ms To prove this lemma it is sufficient to consider five cases of the position of the first letter of the weak square inside the canonical framgment, and for each case easily obtain the contradiction: corresponding letters of the weak square will not match, or increasing\,/\,decreasing on the same word will be different.

\medskip\ms\stat{Lemma (12)}{When morphism $\f$ is applied to a weakly squarefree word, the image of this word does not contain any weak squares of the form $aX\!Xa$, where $|X|$ is divisible by 5.}

\ms The two cases of the position of $X\!X$ inside the sequence of canonical fragments are possible.

\ms The first one: $X$ is an image of several letters --- $X = \f(W)$. We consider the weak square $a\,\f(W)\f(W)\,a$. It means that the first letter $a$ is the last in the canonical fragment, and the second letter $a$ is the first in the canonical fragment. As a consequence, these canonical fragments are both equal to $\f(a)$. Finally, the original word must contain the weak square $aW\!Wa$, which is impossible.

\ms The second one: the words $X$ begin at the same positions inside different canonical fragments. Then these fragments are identical, and the word $X$ must end with the $a$ letter. It means that inside the canonical fragment two successive letters are both equal to $a$, which is impossible.

\medskip\ms Now let the image of a weakly squarefree word contain a weak square $aX\!Xa$. Lemmas 6--12 consider all the cases of the length of $X$, and for each case it is proved that weak squares with this particular length of $X$ do not exist in the image.

\ms The Theorem is proved.


\medskip\ms\stat{Theorem 3}{The morphism $\f: 1 \ra 12321, 2 \ra 23132, 3 \ra 31213$ over the ternary alphabet has a fixed point $\f ^\infty (1)$.}

\ms Indeed, $\f(1) = 1\,2321 = 1W$, and by Lemma (0) the fixed point $\f ^\infty (1)$ exists.

\medskip\ms\stat{Remark}{The morphism $\f: 1 \ra 12321, 2 \ra 23132, 3 \ra 31213$ over the ternary alphabet is not squarefree and not overlap-free.}

\ms Indeed, $\f(212) = 231\ 3\,212\,3\,212\,3\ 132$ contains an overlap, and, consequently, a square.


\ \\
\medskip\ms\stat{Theorem 4}{The morphism $\f: 1 \ra 1221, 2 \ra 2332, 3 \ra 3113$ over the ternary alphabet is a Thue morphism.}

\begin{center}
\includegraphics[width=4.7cm]{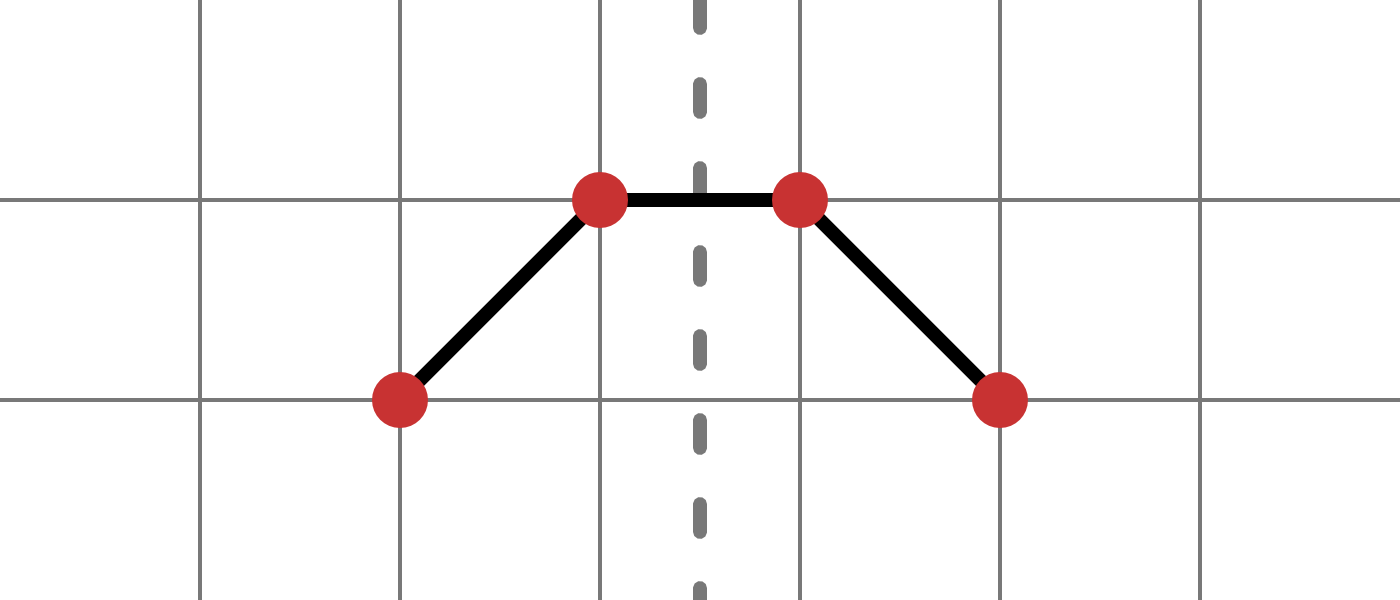}
{\\ \small\it a canonical fragment of the morphism $\f$, \\
increasing and decreasing factors of it}
\end{center}

\ms We need a series of lemmas to prove this theorem.

\ms It is necessary to show that the morphism $\f$ is cubefree, overlap-free and has a fixed point. The morphism $\f$ that we consider has a fixed point ---  $\f ^\infty (1)$ --- because $\f(1) = 1221 = 1W$, and this is exactly what the wording of Lemma (0) requires.

\ms Lemmas (13)--(15) consider a generalized form of cubefreeness and overlap-freeness. They prove that when the morphism $\f$ is applied to a cubefree word the image of it does not contain any overlaps of the corresponding length.

\ms We can consider this form generalized as long as an overlap-free word is always cubefree.

\medskip\ms\stat{Lemma (13)}{When morphism $\f$ is applied to a cubefree word, the image of this word does not contain any cubes of a letter.}

\ms Let the image of the word contain a cube of a letter $aaa$. Then two of these three letters belong to a canonical fragment's edge. On the other hand, letters on the edge of the canonical fragments are different if the morphism $\f$ is considered.

\vfill\eject

\medskip\ms\stat{Lemma (14)}{When morphism $\f$ is applied to a cubefree word, the image of this word does not contain any overlaps of the form $aX\!aX\!a$, where $|X|=1$.}

\ms Let the image of the word contain such an overlap. The length of this overlap is equal to 5. Two successive canonical fragments are sufficient to cover this overlap.

\ms Let us now notice that for any canonical fragment $abcd$ the following equation takes place: \\
\centerline{$b=c=a \op 1=d \op 1$.}

\ms Two $a$ letters out of three in the overlap belong to one canonical fragment. The difference between the position of these letters is equal to 2. However, in the canonical fragment $abcd$ we have $a \not = c$ and $b \not = d$.

\ms As a consequence, we can not have an overlap of length 5 in the image of a cubefree word.

\medskip\ms\stat{Lemma (15)}{When morphism $\f$ is applied to a cubefree word, the image of this word does not contain any overlaps of the form $aX\!aX\!a$, where $|X|=2$.}

\ms Let the image of the word contain such overlap. The length of this overlap is equal to 7. Three successive canonical fragments are sufficient to cover this overlap. Let us denote these fragments by $a_1a_2a_3a_4$, $b_1b_2b_3b_4$ and $c_1c_2c_3c_4$.

\ms Let us now consider the cases of location of two $X$ words inside the three canonical fragments.

\begin{enumerate}
  \item $a_2a_3 = b_1b_2$. These words can not be equal, because $a_2 = a_3$, whereas $b_2 = b_1 \op 1$.
  \item $a_3a_4 = b_2b_3$. These words can not be equal, because $b_2 = b_3$, whereas $a_3 = a_4 \op 1$.
  \item $a_4b_1 = b_3b_4$. Now let us consider the $a$ letters preceding the $X$ words. Therefore, $a = a_3 = a_4 \om 1$, and $a = b_2 = b_3$. On the other hand, these equations can not be true at the same time, because $b_3 = a_4$.
  \item $b_1b_2 = b_4c_1$. Now let us consider the $a$ letters following the $X$ words. Therefore, $a = b_3 = b_2$, and $a = c_2 = c_1 \op 1$. On the other hand, these equations can not be true at the same time, because $b_2 = c_1$.
\end{enumerate}

\ms As a result, there are no overlaps of length 7 in the image of a cubefree word.

\medskip\ms\stat{Lemma (16)}{When morphism $\f$ is applied to a cubefree word, the image of this word does not contain any overlaps of the form $aX\!aX\!a$, where $|aX|$ is divisible by 4.}

\ms Because the morphism $\f$ is cyclic, its canonical fragment can be determined by a letter and its position inside the fragment. Furthermore, it is obvious that the $a$ letters lie on the same positions inside the corresponding canonical fragments. The reason for this is that the global positions of the $a$ letters in the image of the word have the same ramainder of the division by 4.

\ms Furthermore, the corresponding letters of the $X$ words lie on the same positions in the canonical fragments and are equal to each other. Therefore, the corresponding canonical fragments are equal, and there was an overlap in the original word. It contradicts the wording.

\medskip\ms The proof of the fact that the image of a cubefree word does not contain any cubes of the form $X\!X\!X$, where $|X|$ is divisible by $4$, is absolutely the same.

\vfill\eject

\medskip\ms\stat{Lemma (17)}{In a sequence of canonical fragments of the morphism $f$ a factor of length more than 4 and not divisible by 4 can not repeat twice.}

\ms Let the sequence of canonical fragments contain a square $X\!X$. Then each word $X$ contains a border between two canonical fragments next to each other. Furthermore, if we compare these two words $X$, we will see that the borders of canonical fragments are shifted in the second one.

\ms On the other hand, it can be noticed that any shift of canonical fragment of the morphism $\f$ has increasing\,/\,decreasing different from original one. Finally, it means that one factor of $X$ must be increasing and decreasing at the same time, which is impossible.

\ms The Lemma is now proved.

\medskip\ms As a result, the morphism $\f$ is a Thue morphism. If the image of a cubefree word contains a cube, or the image of an overlap-free word contains an overlap, the corresponding lemma proves that there can not be any cubes and overlaps of the corresponding length.

\ms The theorem is proved.

\medskip\ms\stat{Remark}{The morphism $\f: 1 \ra 1221, 2 \ra 2332, 3 \ra 3113$ over the ternary alphabet is not squarefree and not weakly squarefree.}

\ms Indeed, the image of each letter contains a square and is a weak square.


\ \\
\medskip\ms\stat{Theorem 5}{The morphism $\f: 1 \ra 121, 2 \ra 232, 3 \ra 313$ over the ternary alphabet is cubefree and has a fixed point.}

\begin{center}
\includegraphics[width=4.7cm]{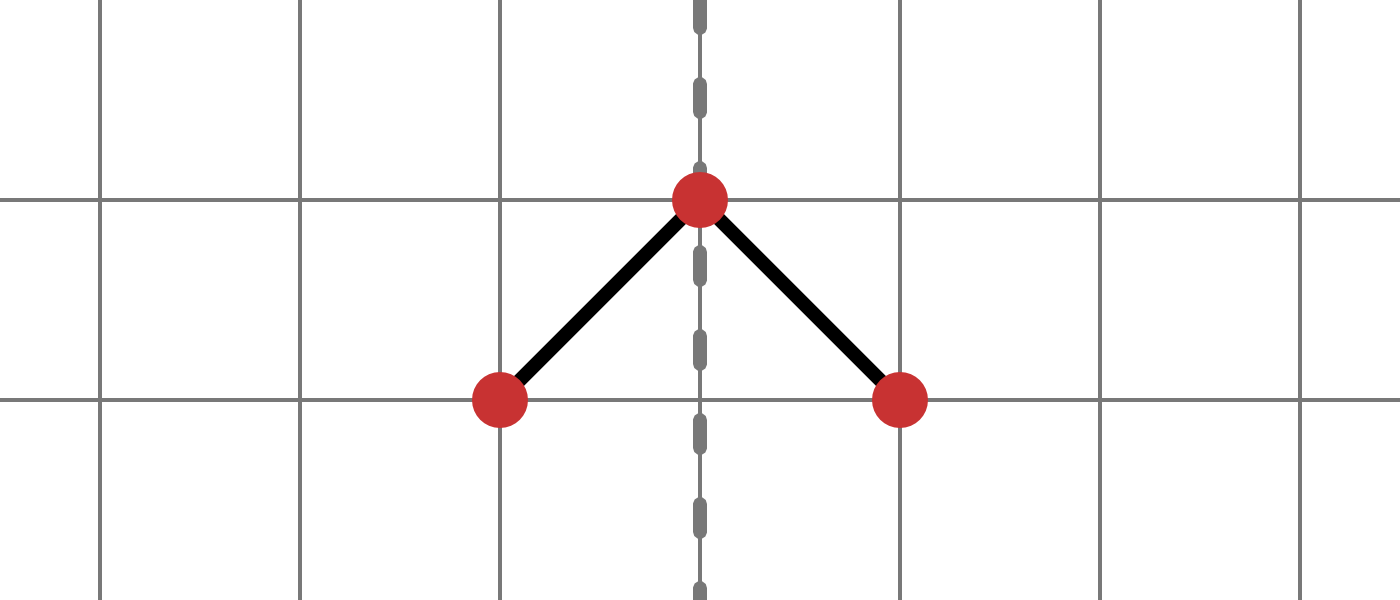}
{\\ \small\it a canonical fragment of the morphism $\f$, \\
increasing and decreasing factors of it}
\end{center}

\ms\stat{Lemma (18)}{When morphism $\f$ is appliaed to a cubefree word, the image of this word does not contain any cubes of letters (of the form $aaa$).}

\ms Let the image of the word contain a cube of a letter $aaa$. At least two out of three letters given belong to a common canonical fragment. However, inside the canonical fragment of the morphism $\f$ successive letters are different.

\medskip\ms\stat{Lemma (19)}{When morphism $\f$ is applied to a cubefree word, the image of this word does not contain any cubes of the form $X\!X\!X$, where $|X|=2$.}

\ms Let the image of the word contain such a cube. Its length is equal to 6. Three successive canonical fragments are sufficient to cover it. Let us denote these canonical fragments by $a_1a_2a_3$, $b_1b_2b_3$ and $c_1c_2c_3$.

\ms Let us consider the cases of the position of the cube inside these canonical fragments.

\begin{enumerate}
  \item $a_1a_2 = a_3b_1 = b_2b_3$. This is impossible because the word $a_1a_2$ is increasing, whereas the word $b_2b_3$ is decreasing.
  \item $a_2a_3 = b_1b_2 = b_3c_1$. This is impossible because the word $b_1b_2$ is increasing, whereas the word $a_2a_3$ is decreasing.
  \item $a_3b_1 = b_2b_3 = c_1c_2$. This is impossible because the word $c_1c_2$ is increasing, whereas the word $b_2b_3$ is decreasing.
\end{enumerate}

\ms As a result, there are no cubes of length 6 in the image of a cubefree word.

\medskip\ms\stat{Lemma (20)}{When morphism $\f$ is applied to a cubefree word, the image of this word does not contain any cubes of the form $X\!X\!X$, where $|X|$ is divisible by 3.}

\ms The proof of this lemma is similar to the proofs of lemmas (5) and (16). Indeed, a canonical fragment can be determined by a letter and its position in the fragment, and the corresponding letters of the words $X$ have the same positions in the corresponding canonical fragments.

\ms As a result, we obtain that the original word contain a cube.

\medskip\ms\stat{Lemma (21)}{When morphism $\f$ is applied to a cubefree word, the image of this word does not contain any cubes of the form $X\!X\!X$, where $|X|$ is greater than 3 and not divisible by 3.}

\ms As long as the canonical fragment of the morphism $\f$ has length 3, and $|X|$ is not divisible by 3, one of the $X$ words begins with the first letter of a canonical fragment, and another one begins with the second letter of a canonical fragment.

\ms This fact can be proved using different methods. For example, it is known that $\{1;2\}$ is the generating set of the group $\mathbb Z / 3 \mathbb Z$.

\ms The fact that one of the $X$ words begins with the first letter of a canonical fragment makes $X\![1]X\![2]$ an increasing word. On the other hand, $X\![1]X\![2]$ must be decreasing, because another $X$ word begins with the second letter of a canonical fragment.

\ms As a result, $X\![1]X\![2]$ must be increasing and decreasing at the same time, which is impossible.

\medskip\ms\stat{Lemma (22)}{The morphism $\f$ has a fixed point $\f ^\infty (1)$.}

\ms We can use lemma (0) to prove it. Indeed, $\f(1) = 121 = 1W$.

\medskip\ms Therefore, the theorem is proved. Lemmas (18)--(21) consider all the cases of the length of the cube and show that the image of a cubefree word can not contain any cubes of such length.

\medskip\ms\stat{Remark}{The morphism $\f$ is neither squarefree nor weakly squarefree nor overlap-free.}

\begin{enumerate}
  \item $\f(12) = 12\,12\,32$ contains a square;
  \item $\f(123) = 12\,1\,23\,23\,1\,3$ contains a weak square;
  \item $\f(212) = 23\,2\,1\,2\,1\,2\,31$ contains an overlap.
\end{enumerate}

\vfill\eject

\section{Low-rank morphisms}

\stat{Lemma (23)}{If $\f$ is a weakly squarefree Thue morphism over the ternary alphabet, then all the first characters of the images of the letters are different; all the last characters of the images of the letters are also different.}

\ms Let us assume that there exist two letters $a$, $b$ such as \\
\centerline{$\f(a) = xU$, $\f(b)=xW$.}

\ms Then $\f(aab)$ contains an overlap $xUxUx$, which contradicts the wording.

\ms Similarly, if $\f(a) = Ux$, $\f(b)=Wx$, then $\f(abb)$ contains an overlap $xWxWx$.

\ms Let us now notice that we appealed only to the fact that a Thue morphism is always overlap-free. Thus, it can be said that {\it if $\f$ is an overlap-free morphism over the ternary alphabet, then all the first characters of the images of the letters are different; all the last characters of the images of the letters are also different.}

\medskip\ms\stat{Lemma (24)}{If $\f$ is a weakly squarefree Thue morphism over the ternary alphabet, then for each $a \in \A$ the first and the second letters of $\f(a)$ are different, the last and the penultimate letters of $\f(a)$ are also different.}

\ms Lemma (23) implies that for every $a \in \A = \{1,2,3\}$ there is a canonical fragment beginning with it. Therefore, if $\f(a)=Uxx$, then there exists $b$ such as $\f(b)=xW$. As a result, $\f(ab)$ contains a cube, which contradicts the wording.

\ms Another proposition of the lemma can be proved similarly.

\medskip\ms\stat{Lemma (25)}{If $\f$ is a weakly squarefree Thue morphism over the ternary alphabet, then for each $a \in \A$ the first character of $\f(a)$ is equal to the last character of $\f(a)$.}

\ms Let there be a letter $a$ such as $\f(a)=xWy$. Then there exists a letter $b$ such as $\f(b)=yU$, and $b \neq a$. As a result, $\f(ab) = xWyyU$ contains a weak square $yy$.

\medskip\ms\stat{Theorem 6}{Over the ternary alphabet there are no uniform weakly squarefree Thue morphisms with a rank lower than 5.}

\ms\stat{Proposition 1}{Over the ternary alphabet there are no uniform overlap-free morphisms with rank 2.}

\ms Let there be an overlap-free morphism $\f: a_1 \ra ab, a_2 \ra cd, a_3 \ra ef$. Lemma (23) implies that $a \neq c \neq e$ and $b \neq d \neq f$.

\ms Let us now consider the image of an overlap-free word $a_1a_1$: $\f(a_1a_1)=abab$. In order that $aba$ is not an overlap, $b \neq a$ is necessary. Similarly it can be proved that $d \neq c$ and $f \neq e$.

\ms Only two types of morphisms satisfy the listed properties: \\
\centerline{$\f_1: a_1 \ra 12, a_2 \ra 23, a_3 \ra 31$;}
\centerline{$\f_2: a_1 \ra 13, a_2 \ra 21, a_3 \ra 32$.}

\ms Let us now show that $\f_1$ and $\f_2$ are not overlap-free.
\begin{enumerate}
  \item $\f_1(a_1a_3a_2a_1a_3a_2) = 123123123123$ contains a cube and, consequently, an overlap.
  \item $\f_2(a_1a_2a_3a_1a_2a_3) = 132132132132$ contains a cube and, consequently, an overlap.
\end{enumerate}

\vfill\eject


\medskip\ms\stat{Proposition 2}{Over the ternary alphabet there are no uniform overlap-free morphisms with rank 3.}

\ms Let there be a weakly squarefree Thue morphism $\f$ with rank 3. The morphism $\f$ has a fixed point, therefore, image of one of the letters begins with this letter. For example, let $\f(1)$ begin with 1. Lemma (24) implies that $\f(1)=12u$ or $\f(1)=13u$; these cases are equivalent, let us consider the first one.

\ms Lemma (25) implies that $\f(1)=121$. Therefore, the canonical fragment that begins with 2 can not have 2 nor 1 in its midde. As a consequence, $\f(a_2)=232$. Similarly, $\f(a_3)=313$. As a result,
\begin{center}
  $\f(1) = 121$ \\
  $\f(a_2) = 232$ \\
  $\f(a_3) = 313$.
\end{center}

\ms However, $\f(a_31a_2) = 31\,3\,12\,12\,3\,2$ contains a weak square.

\ms Finally, the proposition is proved.

\medskip\ms\stat{Proposition 3}{Over the ternary alphabet there are no uniform overlap-free morphisms with rank 4.}

\ms Let there be a weakly squarefree Thue morphism $\f$ with rank 4. Similarly to the previous proposition, only the case $\f(1) = 12mn$ may be considered. This time, by lemma (25), $\f(1)=1231$.

\ms It is easy to see that $\f(a_2)$ is equal to 2132 or to 2312; $\f(a_3)$ is equal to 3123 or to 3213.

\ms The case $\f(a_2) = 2312$ is impossible: $\f(1a_2) = 1\,2\,31\,2\,31\,2$ contains an overlap. Therefore, $\f(a_2) = 2132$. Let us notice that now $\f(a_3) = 3213$ is impossible, because $\f(a_3a_2) = 3\,21\,3\,21\,3\,2$ contains an overlap.

\ms As a result, we have the only possible morphism:
\begin{center}
  $\f(1) = 1231$ \\
  $\f(a_2) = 2132$ \\
  $\f(a_3) = 3123$.
\end{center}

\ms However, in this case $\f(a_31) = 3\,1\,23\,1\,23\,1$ contains an overlap.

\medskip\ms As a result, the cases of ranks $2$, $3$, $4$ are considered, and the theorem is proved.

\medskip\ms\stat{Lemma (26)}{Over the binary alphabet there are no cubefree and weakly squarefree words of length greater than 5.}

\ms Indeed, if the word is weakly squarefree, that it does not contain any squares of letters. It means that after the letter 1 it is always 0, and after 0 it is always 1.

\ms As a result, when the length is equal to 6, we will have a cube 101010 or 010101.

\medskip\ms\stat{Theorem 7}{There are no weakly squarefree Thue morphisms over the binary alphabet.}

\ms Let there be a weakly squarefree Thue morphism $\f$ over the binary alphabet. It has a fixed point $A = \f ^\infty (a)$. Let us denote the sequence that generates this fixed point by $A_n$.

\ms The limit of $|A_n|$ is equal to infinity. This implies that at some point $|A_N|$ will be greater than 6.

\ms The word $a$ is weakly squarefree and overlap-free. However, a word $A_N$  contains a weak square or an overlap. Therefore, $\f$ is not weakly squarefree or cannot be viewed as a Thue morphism.

\vfill\eject


\section{Weak squarefreeness}

\newstr\stat{Theorem 8}{A cyclic squarefree morphism over the ternary alphabet is always weakly squarefree.}

\ms To prove this theorem we need a number of lemmas.

\medskip\ms\stat{Lemma (27)}{If a squarefree morphism is applied to a word and if the image of this word contains a weak square, then the original word contains a square.}

\ms Indeed, a word containing a weak square contains a square. However, if squarefree morphism is applied to a squarefree word, then the image must be squarefree, but it is not.

\medskip\ms\stat{Lemma (28)}{If the morphism $\f$ is cyclic and squarefree, then for each letter $a \in \A$ the first and the last letters of $\f(a)$ are equal.}

\ms Let there be $a \in \A$ such that its image begins with the letter $x$ and ends with the letter $y$. Arising from the properties of the ternary alphabet, $y$ is equal to $x \op 1$ or $x \om 1$.

\ms Let us consider the case when $y = x \op 1$. As long as the morphism $\f$ is cyclic, the image of $a \op 1$ begins with the letter $y$. Therefore, the image of a squarefree word $\f(a(a \op 1))$ contains a square.

\ms On the other hand, $\f$ is squarefree.

\medskip\ms\stat{Lemma (29)}{If the morphism $\f$ is squarefree and cyclic, then the image of a weakly squarefree word cannot contain a square $W\!W$ such as $|W\!|$ is not divisible by the rank of the morphism $\f$.}

\ms Let us denote the rank of the morphism $\f$ by $L$.

\ms Let the image of a weakly squarefree word contain such square. Then $|W\!|$ can be greater or lower than $L$. Let us consider these cases separately.

\begin{enumerate}

\item If $|W\!|<L$, then the square $W\!W\!$ can be covered with three successive canonical fragments. Indeed, if this square cannot be covered with three canonical fragments, then its length must exceed the length of two canonical fragments, which is equal to $2L$. As a consequence, we have $2|W\!|>2L$ $\Longrightarrow$ $|W\!|>L$. Hovewer, we consider the case when $|W\!|<L$.

\ms Let us now notice that if a word has the length lower or equal to 3, then weak squarefreeness implies squarefreeness. As long as the morphism $\f$ is squarefree, the image of a weakly squarefree word of length 3 can not contain the square $W\!W\!$, which contradicts the assumption.

\item Let now $|W\!|$ be greater than $L$. We are going to show that the increasing and decreasing on the canonical fragment is equivalent to the increasing and decreasing on a shift of the canonical fragment. We need to prove that:

\qquad\qquad $\exists h:\ \ \forall\, i = \ovl{1,L-1}\ \ \forall\, a \in A$

\smallskip\qquad\qquad\qquad $\bigl(\f(a)\bigr)[i+1] = \bigl(\f(a)\bigr)[i] \op 1\ 
\ \Longrightarrow$

\qquad\qquad\qquad $\Longrightarrow\ \ \bigl(\f(a)\bigr)[(i+h+1)\ \mathrm{mod}\ L] 
=
  \bigl(\f(a)\bigr)[(i+h)\ \mathrm{mod}\ L] \op 1$;

\smallskip\qquad\qquad\qquad and, similarly,

\qquad\qquad\qquad $\bigl(\f(a)\bigr)[i+1] = \bigl(\f(a)\bigr)[i] \om 1\ \ 
\Longrightarrow$

\qquad\qquad\qquad $\Longrightarrow\ \ \bigl(\f(a)\bigr)[(i+h+1)\ \mathrm{mod}\ L] 
=
  \bigl(\f(a)\bigr)[(i+h)\ \mathrm{mod}\ L] \om 1$.

\medskip\centerline{\small\it (in other words, in the canonical fragment the relation between two letters}
\centerline{\small\it is equivalent to the relation between the letters that are shifted by $h$ positions)}
\smallskip\centerline{\small\it (by $U[p]$ we denote the $p$--th letter of the word $U$)}

\vfill\eject

\ms To prove lemma it is sufficient to provide the $h$ number. Let $h$ be equal to $(|W\!|\ \mathrm{mod}\ L)$.

\ms Let us denote the square $W\!W$ by $W_1W_2$,\ $W_1 = W_2 = W\!$. And let us consider the relation between $i$--th and $i+1$--th letters of the canonical fragment. Let \\
\centerline{$\bigl(\f(a_j)\bigr)[i+1] = \bigl(\f(a_j)\bigr)[i] \op 1$.}

\ms As long as $|W_1|>L$, for each $i$ it is true that $i$--th and $i+1$--th letters of some canonical fragment will be successive in $W_1$ or $W_2$. In other words,

\centerline{$\exists\, j, k :\ \ W_1[k] \equiv \bigl(\f(a_j)\bigr)[i],
  \ \ W_1[k+1] \equiv \bigl(\f(a_j)\bigr)[i+1]$;}
\smallskip\centerline{$W_1[k+1] = W_1[k] \op 1$.}

\ms Furthermore, if $k$--th and $k+1$--th letters of the word $W_2$ belong to one canonical fragment, then they have the same relation: $W_2[k+1] = W_2[k] \op 1$. However, the positions of $W_2[k]$ and $W_2[k]$ in the canonical fragment will be different from $i$ and $i+1$.

\ms Finally, this difference will be equal to the remainder of the division $|W\!|$ by $L$ --- which is equal to $h$. We have:

\centerline{$\bigl(\f(a_m)\bigr)[(i+h+1)\ \mathrm{mod}\ L] =
  \bigl(\f(a_m)\bigr)[(i+h)\ \mathrm{mod}\ L] \op 1$.}
\smallskip\centerline{\small\it(we need the eternal $\mathrm{mod}\ L$ in case the numbers}
\centerline{\small\it $i$ or $i+h$ are greater than $L$)}

\ms Now let us notice that the relations of the form $\op 1\ /\ \om 1$ on the same positions take place in all the canonical fragments of the cyclic morphism. Therefore, if what we need is true in the canonical fragments $\f(a_j)$ or $\f(a_m)$, then it is true in the canonical fragment $\f(a)$ for any $a$.

\ms The preservation of the increasing\,/\,decreasing for the shift is shown.

\begin{center}\includegraphics[width=7.3cm]{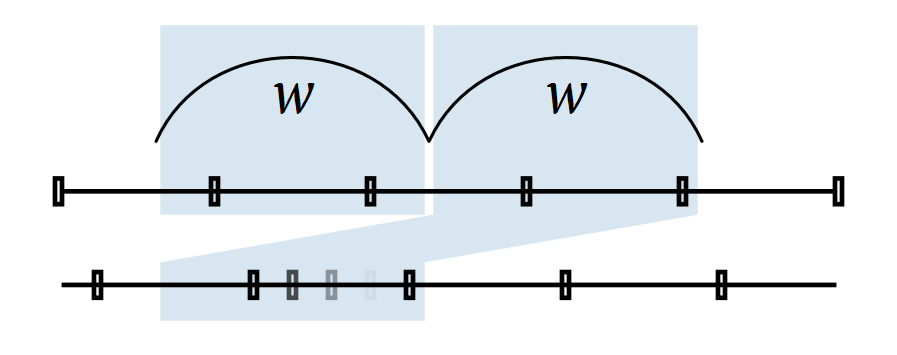}\end{center}

\ms Now we know that the relation between $i$--th letter of the canonical fragment and the next one is the same for the $(i+h)$--th letter of the canonical fragment and the next one. So, it is the same for $(i+2h)$, $(i+3h)$, ..., $(i+kh)$--th letters of the canonical fragment and the next ones.

\ms In other words, it is sufficient to define the segment of the first $h'$ letters of the canonical fragment to know that on other segments of this length the increasing\,/\,decreasing is similar to the increasing\,/\,decreasing on the first one.

\ms From the number theory, $h'$ is equal to the greatest common divisor of $h$ and $L$. We know that $h<L$, so $h'$ is the divisor of $L$ lower than $L$.

\ms Now the canonical fragment is divided into segments of length $h'$, and monotony on all this segments is similar. Furthermore, the next segment is always connected to the previous one with one relation $\op 1\ /\ \om 1$, which is also similar everywhere. It is easy to determine this relation by the first, $h'$--th and $h'+1$--th letters of the canonical fragment.

\vfill\eject

\ms Let us only notice that two successive segments can not be equal, because in this case the image of a letter will contain a square.

\ms Now let us consider the canonical fragment $\f(a)$.

\begin{itemize}

    \item[(a)] Let the canonical fragment be divided into two segments. Let them be equal to $S$ and $S \op 1$. Then the images of two other letters have forms $f(a') = (S \op 1)(S \om 1)$ and $\f(a'') = (S \om 1) S$. In this case $\f(a'aa'')=(S \op 1)(S \om 1)S\,(S \op 1)(S \om 1)S$ contains a square. On the other hand, the morphism $\f$ is squarefree.

    \item[(b)] Let us now have more than two segments --- three, at least. Furthermore, let the next segment always be equal to the previous one increased by 1.

\ms Let us now consider a canonical fragment, and let us now take another fragment such that the sequence of segments in it continues the sequence of segments in the first one. It will be another canonical fragment, because the last letter of the first canonical fragment and the first letter of the second canonical fragment will be different, and we have the lemma (28).

\ms Then, on the border of two canonical fragments we have \\
\centerline{$...\ S(S \op 1)(S \om 1) \ S(S \op 1)(S \om 1)\ ...$,}
which is a square. On the other hand, the morphism $\f$ is squarefree.

\end{itemize}
\end{enumerate}

\ms The lemma is now proved.

\medskip\ms\stat{Lemma (30)}{If the morphism $\f$ is squarefree and cyclic, then the image of any word does not contain any weak squares of the form $aX\!Xa$, where $|X|$ is divisible by $L$, and the borders of $X$ are not borders of any canonical fragment.}

\ms We already know that the canonical fragment of a cyclic morphism can be determined by a letter and its position in the canonical fragment. So, if two edges and the border between the $X$ words belong to the canonical fragments $\f(a_i)$, $\f(a_j)$ and $\f(a_k)$, then these borders lie on the same positions, and the canonical fragments are similar.

\ms It means that the $X$ word ends with the $a$ letter, because the end of $X$ and $a$ have the same positions in similar canonical fragments. Therefore, $X$ also begins with the $a$ letter.

\ms As a result, there is a square $aa$ inside the canonical fragments, which contradicts squarefreeness of $\f$. The lemma is proved.

\medskip\ms Finally, if $\f$ is a cyclic squarefree morphism, than the square in the image of a weakly squarefree word can only be the image of a square. Let us now consider two canonical fragments next to the square. If the square is weak, then the last letter of the first of these fragment is equal to the first letter of the second one.

\ms So, these fragments are similar. Consequently, there was a weak square in the original word, which contradicts the wording.

\medskip\ms\stat{Corollary}{As long as the Leech morphism is cyclic and squarefree, it is also weakly squarefree.}

\ms As a result, the Leech morphism has all the considered properties: it is cubefree, overlap-free, squarefree and weakly squarefree. It is also cyclic and has a fixed point.

\medskip\ms\stat{Remark}{There is a weakly squarefree morphism over the binary alphabet:}\\
\centerline{$\f: 0 \ra 01, 1 \ra 01$.}

\ms The norphism $\f: 0 \ra 01, 1 \ra 10$ is not weakly squarefree: $\f(10)$ contains a weak square. The theorem 7 proves that there are no weakly squarefree Thue morphisms over the binary alphabet.

\ms The theorem is now proved.


\ \\
\section{Leech's morphism and optimal ranks}

\newstr\stat{Theorem 9}{Over the ternary alphabet there are exactly 144 unifrom squarefree morphisms of rank 11; there are no uniform squarefree morphisms with a lower rank.}

\ms This result was obtained using a program written in the Objective Caml language. In the paper [4] there is a criterion of squarefreeness of an uniform morphism over the ternary alphabet. The criterion says that for a morphism to be squarefree it is sufficient to have squarefree images of the squarefree words of length 3.

\ms Squarefreeness of a word is examined using the function {\tt squarefree\_check}. The computational complexity of it is equal to $O(n^2)$, but it is insignificant for us.

\begin{verbatim}
let squarefree_check str =
  let res = [|true|] in
  let l = String.length str in
    for i=0 to l-2 do (
      for m=1 to (l-i)/2 do (
      if (String.sub str i m = String.sub str (i+m) m)
        then (res.(0) <- false)
      ) done
    ) done;
  res.(0);;
\end{verbatim}

\ms The array {\tt res} contains the intermediate result, and the line

{\tt \quad\qquad if (String.sub str i m = String.sub str (i+m) m)

\quad\qquad then (res.(0) <- false)}

\newstr checks if there is a square of length $2m$ beginning at the $i$--th position.

\ms The array of all the squarefree words of length 11 was created, and the images of letters were take nfrom this array. Then the squarefreeness check was made, and if the morphism happened to be squarefree, the images of letters were printed on the screen.

\ms The source code can be seen in the Appendix.

\medskip\ms\stat{Proposition}{Over the ternary alphabet there are exactly 2 unifrom squarefree morphisms of rank 11 such that all the other squarefree morphisms of rank 11 can be obtained from them using the natural operations. Namely:}\\
\centerline{$\f_1: 1 \ra 12131232123, 2 \ra 13212321323, 3 \ra 13213121323$;}
\centerline{$\f_2: 1 \ra 12313231213, 2 \ra 12321231213, 3 \ra 23132123213$.}

\ms Let us prove this fact. We denote a arbitrary permutation of the set $\{1,2,3\}$ by $\sigma$, and let $\ovl{\sigma}$ be the morphism induced by this permutation. Furthermore, let the inversion --- the function that returns the word in reverse --- be denoted by $\iota$.

\ms If the morphism $\f$ is uniform and squarefree, then the morphisms $\iota \circ \f$, $\f \circ \ovl{\sigma}$ and $\ovl{\sigma} \circ \f$ are also uniform and squarefree.

\ms For a particular morphism $\f$ there are exactly $6 \cdot 6 = 36$ morphisms of the form $\ovl{\sigma_1} \circ \f \circ \ovl{\sigma_2}$ and 36 morphisms of the form $\iota \circ \ovl{\sigma_1} \circ \f \circ \ovl{\sigma_2}$.

\ms It is shown manually that for two uniform squarefree morphisms of rank 11: \\
\centerline{$\f_1: 1 \ra 12131232123, 2 \ra 13212321323, 3 \ra 13213121323$;}
\centerline{$\f_2: 1 \ra 12313231213, 2 \ra 12321231213, 3 \ra 23132123213$,}
\newstr all the other uniform squarefree morphisms of rank 11 have the form $\ovl{\sigma_1} \circ \f_i \circ \ovl{\sigma_2}$ or $\iota \circ \ovl{\sigma_1} \circ \f_i \circ \ovl{\sigma_2}$.

\medskip\ms\stat{Remark}{The lower bound for the rank of the squarefree morphisms over the ternary alphabet is also provided in the article [8].}

\ms In the article [8] a morphism obtained by inversion of $\f_2$, is also provided. We have to remark that in this article the result is obtained independently. Here we also list all the morphisms that reach the lower bound of the rank.

\medskip\ms\stat{Theorem 10}{Over the ternary alphabet there are no cyclic squarefree morphisms of a rank lower than 13.}

\ms The Proposition shows that an uniform squarefree morphism of rank 11 over the ternary alphabet is not cyclic. For the morphisms of rank 12 we can also list all the squarefree ones and check their cyclicness.

\medskip\ms\stat{Theorem 11}{Over the ternary alphabet there are no morphisms of a rank lower than 13 that are squarefree, cubefree and overlap-free at the same time.}

\ms In the aricles [6]--[8] the criteria of overlap- and cubefreeness of a morphism over the ternary alphabet are provided. The author created the programs similar to those for squarefreeness that prove the given fact.


\ \\
\section{Letter combinations in squarefree words}

\newstr Let {\it two-letter combination} be an arbitrary word of length 2. Over the ternary alphabet there are 6 squarefree two-letter combinations: $ab$, $ac$, $ba$, $bc$, $ca$ and $cb$. Let us try to answer the following question: can we do without one of squarefree two-letter combinations while constructing a squarefree word?

\ms\stat{Theorem 12}{Over the ternary alphabet a squarefree word of length more than 13 contains all the squarefree two-letter combinations.}

\ms Let us prove the theorem. We will try to construct a squarefree word without the two-letter combination $ab$.

\ms Let us notice that it is impossible to construct a squarefree word of length more than three using only letters $b$ and $c$. It means that not later than in the 4--th position there will be an $a$ letter. We know that after the $a$ letter there can be only $c$, and after $b$ or $c$ --- any letter but the given one. Guided by this rules, let us construct all the squarefree words until we obtain words with squares.

\vfill\eject
\begin{center}
\includegraphics[width=10.2cm]{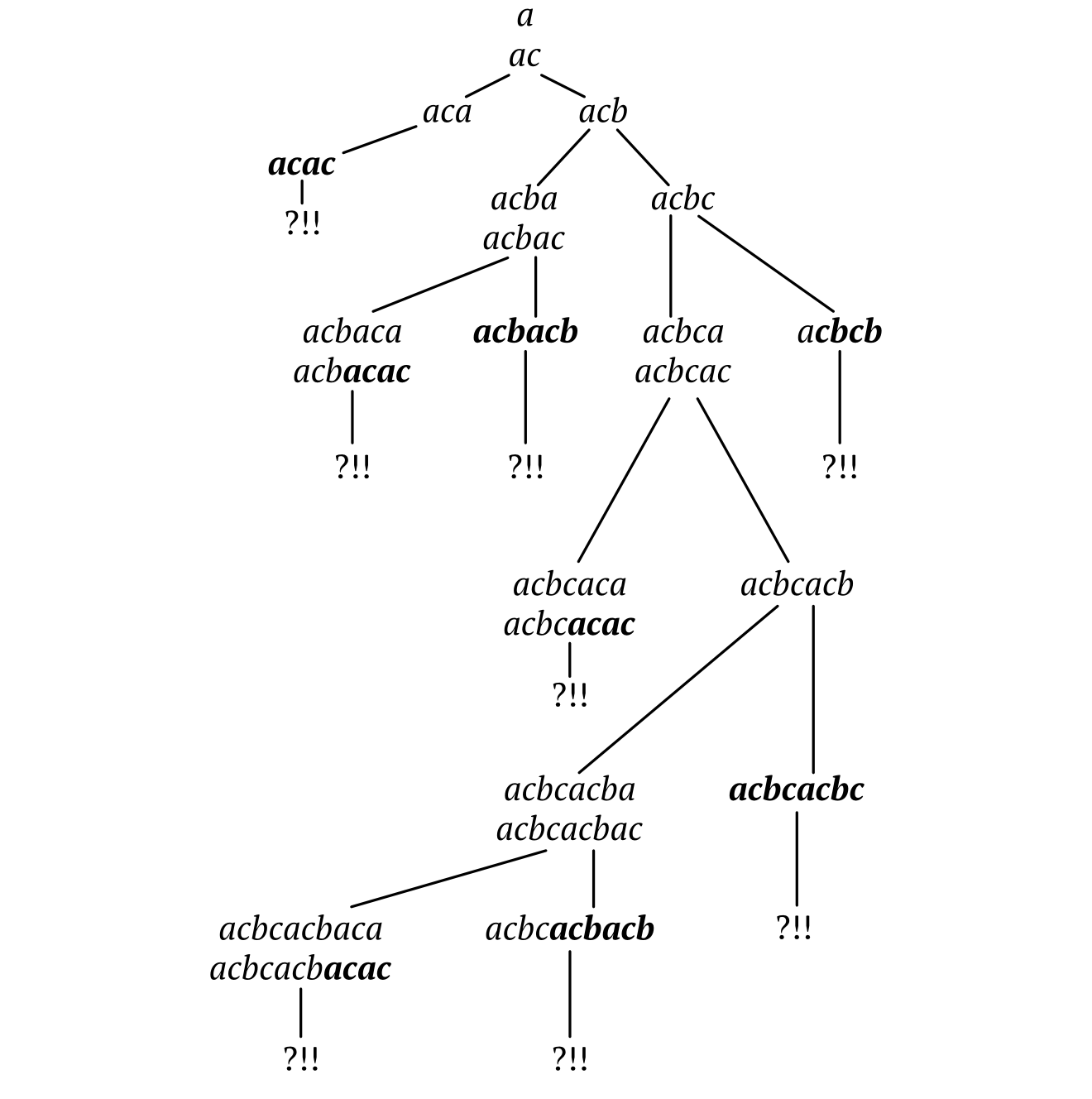}
\end{center}

\ms As a result, $bcba\,cbca\,cbaca$ is the longest squarefree word without the combination $ab$. Its length is equal to thirteen.

\ms The theorem is proved.

\medskip\ms\stat{Remark}{Over a more than three - letter alphabet there are squarefree words of any length without an arbitrary two-letter combination.}

\ms Let us show that there exist squarefree words of any length without the combination $ab$.

\ms The given alphabet without the letter $a$ is not less than ternary. Therefore, we can construct squarefree words of any length over it. These words will not contain the $a$ letter, and, consequently, all the combinations with it.

\ms The remark is proved.

\medskip\ms Similarly to the two-letter combinations we will define the three-letter combinations. A three-letter combination is an arbitrary word of lenth 3. Let us try to answer the following question: which three-letter combinations are necessary to exist inside squarefree words?

\medskip\ms\stat{Theorem 13}{Over the ternary alphabet a squarefree word of length more than 36 contains all the squarefree three-letter combinations of different letters (of the form $abc$).}

\ms Let us again try to construct a squarefree word without the combination $abc$. The theorem 12 implies that not later than in the 12--th position there will be the two--letter combination $ab$. We know that after it there can be only the letter $a$. Now let us continue constructing a squarefree word.

\vfill\eject
\begin{center}\includegraphics[width=13cm]{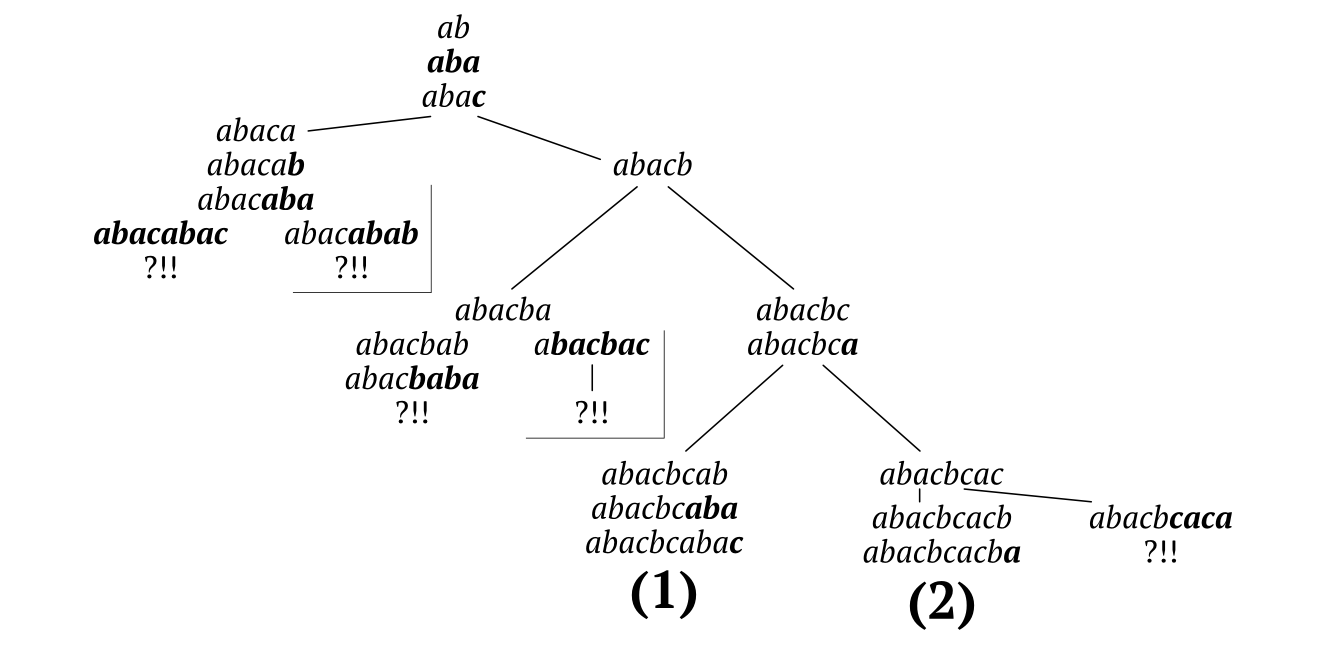}\end{center}

\ms As a result, now we know that a squarefree word beginning with $ab$ and not containing the combination $abc$ can only have forms $abacbcabac..$ or $abacbcacba..$\ . Let us consider these cases separately.

\ms In the first case we always obtain a square:
\begin{center}\includegraphics[width=14.2cm]{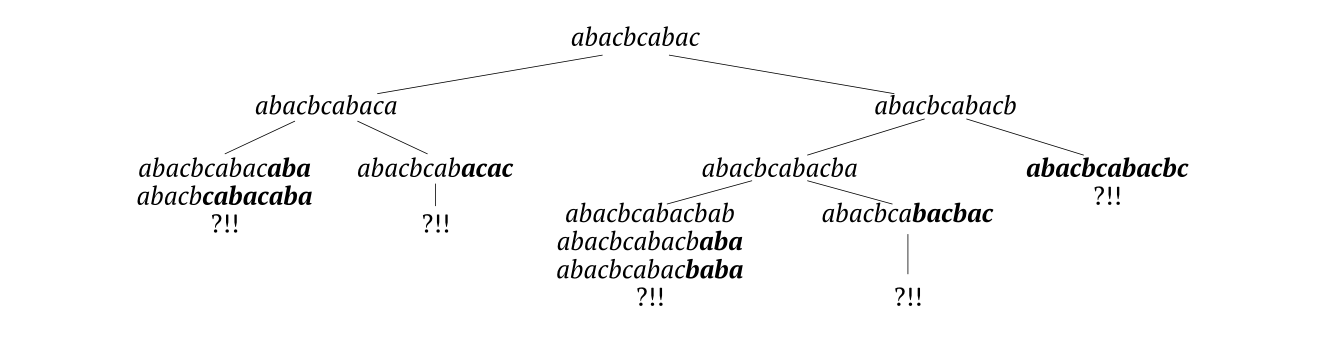}\end{center}

\ms Now the second case is left.
\begin{center}\includegraphics[width=14.2cm]{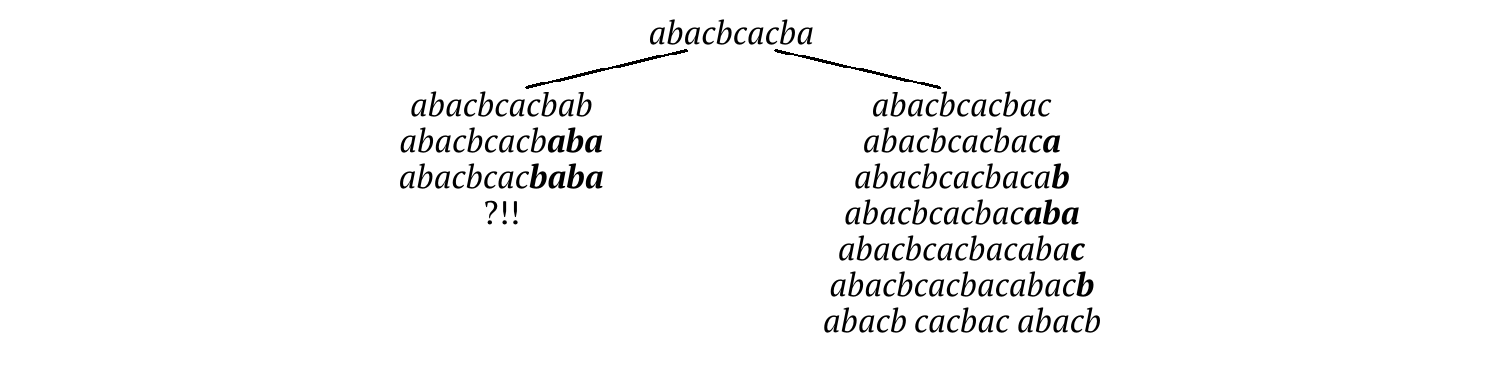}\end{center}

\ms Let us notice that the only long squarefree word ends with $abacb$. Going back to the main tree of cases, we can see that from $abacb$ we always obtain words with squares or the word we are considering.

\ms The longest squarefree word obtained from $abacb$ is $aba\,cbc\,aba\,cbab$ with its length being equal to 13. Adding 13 to the length of the word $abacb\,cacbac$ we get 24. Finally, before the first $ab$ there can be up to 12 letters.

\ms As a result, we have the upper bound, which is equal to 36.

\medskip\ms\stat{Proposition}{Over the ternary alphabet there are squarefree words of any length without the three-letter combination $aba$.}

\ms Let us provide an example of an infinite squarefree word without the factor $aba$. The known Thue result [2] which is the infinite word generated by the morphism
  $$\f: a \ra abcab, b \ra acabcb, c \ra acbcacb$$
does not contain the three-letter combination $cbc$. Indeed, this factor cannot appear inside the images of the letters nor on their borders.

\ms Furthermore, we provide a squarefree word of length 718 without two three-letter combinations: $aba$ and $bab$:

\begin{center}
{\footnotesize\tt

\begin{tabular}{|c|}
\hline
abcbacbcacbacabcbacbcacbabcbacbcabcbacabcacbcabcbacbcacbacab \\
cbacbcabcbacabcacbcabcbacbcacbabcbacbcabcbacabcacbcabcbabcac \\
bacabcbacbcacbacabcbabcacbcabcbacbcacbacabcbacbcabcbacabcacb \\
cabcbacbcacbacabcbacbcabcbabcacbcabcbacbcacbacabcbabcacbcabc \\
bacabcacbcabcbabcacbacabcbacbcacbacabcbabcacbcabcbacbcacbaca \\
bcbacbcabcbacabcacbcabcbacbcacbacabcbacbcabcbabcacbcabcbacbc \\
acbacabcacbcabcbacbcacbabcbacbcabcbacabcacbcabcbacbcacbacabc \\
bacbcabcbacabcacbcabcbacbcacbabcbacbcabcbacabcacbcabcbabcacb \\
acabcbacbcacbacabcbabcacbcabcbacbcacbacabcbacbcabcbacabcacbc \\
abcbacbcacbacabcbacbcabcbabcacbcabcbacbcacbacabcbabcacbcabcb \\
acabcacbcabcbabcacbacabcbacbcacbacabcbabcacbcabcbacbcacbacab \\
cbacbcabcbacabcacbcabcbacbcacbacabcbacbcabcbabcacbcabcbaca \\
\hline
\end{tabular}
}
\end{center}

\medskip\ms\stat{Remark}{Over a more than three - letter alphabet there are squarefree words of any length without an arbitrary three-letter combination.}

\ms Indeed, we can construct a squarefree word without the combination $abc$ or $aba$ just constructing a squarefree word over a not less than ternary subalphabet without the $a$ letter.

\medskip\ms\stat{Remark}{Finiteness of the trees of cases can be checked using the computer, but the manual check is more of interest.}


\ \\
\section{The applications of theory and the results}

The theory of non-repeating words found its applications in many branches of mathematics and other sciences. In the XX century mathematicians Novikov and Adyan solved the Burnside problem of local finiteness of the variety of periodical groups.

\ms Furthermore, Max Euwe, who is a mathematician and an ex -- world champion in chess, rediscovered the Thue-Morse sequence and provided the example of an infinite chess game that cannot be flattened draw. This led to a change to the rules.

\ms Non-repeating sequences can be applied in the game theory, the number theory, cryptography, symbol dynamics (the theory of billiards, for example) and even in genetics. The evidence that combinatorics on words is important is that the Thue-Morse sequence had been rediscovered several times before the Thue works were publishid in Europe. Many articles are published in this field these days as well.

\ms In this work the questions we addressed have been answered. New sets of non-repeating words and morphisms are obtained. Also, the optimal evaluations of the ranks of the morphisms were provided.


\vfill\eject
\section{References}

\def\sqb#1{[#1]}

\begin{itemize}
  \item[\sqb{1}] Axel Thue, Uber unendliche Zeichenreihen; Norske Vid. Selsk Skr. I Mat.-Nat. Kl.; Christiania; 1906 --- 1--22.
  \item[\sqb{2}] Axel Thue, Uber die gegenseitige Lage gleicher Teile gewisser Zeichenreihen; Norske Vid. Skrifter I Mat.-Nat. Kl.; Christiania; 1912 --- 1--67.
  \item[\sqb{3}] John Leech, A problem on strings of beads; Math. Gazette 41; 1957 --- 277--278.
  \item[\sqb{4}] Max Сrochemore, Sharp characterizations of squarefree morphisms; Theor. Comput. Sci.; 18 \linebreak (1982) --- 221--226.
  \item[\sqb{5}] Jean Berstel, Some recent results on squarefree words; Lecture Notes in Computer Science; 166(1984) --- 14--25.
  \item[\sqb{6}] G. Richomme, F. Wlazinski, Existance of finite test-sets for k-power-freeness of uniform morphisms --- arxiv.org/pdf/cs/0512051.
  \item[\sqb{7}] Christopher R. Tompkins, Overlap-free Morphisms; University of South Carolina; Mathematics; 2009 --- 1--67.
  \item[\sqb{8}] F. Brandenburg, Uniformly growing k-th powerfree homomorphisms; Theor. Comput. Sci.; 23(1983) --- 69--82.
\end{itemize}

\vfill\eject


\section{Appendix (A). The source code}

Let us provide the source code of the program, obtaining all the squarefree morphisms of rank 11.

{\small
\begin{verbatim}
(* Objective Caml v 4.00.1 *)
let int_of_bool b = if b=true then 1 else 0;;

let squarefree_check str =
  let res = [|true|] in
  let l = String.length str in
    for i=0 to l-2 do (
      for m=1 to (l-i)/2 do (
      if (String.sub str i m = String.sub str (i+m) m)
        then (res.(0) <- false)
      ) done
    ) done;
  res.(0);;

let printl lis = List.iter (fun x -> print_string (x^" ")) lis;;
let extend str = [str^"1"; str^"2"; str^"3"];;
let step lis = List.flatten (List.map extend lis);;
let rec nextend n = if n=0 then [""] else step (nextend (n-1));;
  (* Generating the list of all the words of length n *)

let main = Array.of_list (List.filter squarefree_check (nextend 11));;
  (* The array of all the squarefree words of length 11 *)

let morphism_application str q w e =
  let apply_morphism t = match t with
    '1' -> main.(q) |
    '2' -> main.(w) |
    '3' -> main.(e) |
    q -> "" in
  let to_flatten = List.map apply_morphism str in
  String.concat "" to_flatten;;
  (* Function — applying a morphism to a word; the images of the letters are taken from the array *)

let ternary = [['1';'2';'1'];['1';'2';'3'];['1';'3';'1'];['1';'3';'2'];
  ['2';'1';'2'];['2';'1';'3'];['2';'3';'2'];['2';'3';'1'];
  ['3';'1';'3'];['3';'1';'2'];['3';'2';'3'];['3';'2';'1']];;
  (* All the squarefree words of length 3 over the ternary alphabet *)

let morphism_check q w e =
  let after_apply = List.filter squarefree_check
    (List.map (fun x -> morphism_application x q w e) ternary) in
  if (List.length after_apply) = 12
    then Printf.printf "%s %s %s \n" main.(q) main.(w) main.(e)
    else ();;
  (* The check of the squarefreeness of the morphism *)

let l = Array.length main;;
for q=0 to l-1 do (
  for w=0 to l-1 do (
    for e=0 to l-1 do (
      morphism_check q w e;
  ) done  ) done ) done;;
\end{verbatim}}

\vfill\eject

\section{Appendix (B). Squarefree morphisms of Rank 11}
Let us provide the list of squarefree morphisms of rank 11 obtained using the written program.
\medskip \\

\centerline{\small\tt
\begin{tabular}{|c|c|c|c|c|c|c|}
\hline
$\f(1)$ & $\f(2)$ & $\f(3)$ & & $\f(1)$ & $\f(2)$ & $\f(3)$ \\
\hline
12131232123 & 13212321323 & 13213121323 & & 12131232123 & 13213121323 & 13212321323\\
12132123213 & 12132313213 & 12321231323 & & 12132123213 & 12321231323 & 12132313213\\
12132313213 & 12132123213 & 12321231323 & & 12132313213 & 12321231323 & 12132123213\\
12312131232 & 12313231232 & 13121323132 & & 12312131232 & 13121323132 & 12313231232\\
12313231213 & 12321231213 & 23132123213 & & 12313231213 & 12321312132 & 12321323132\\
12313231213 & 12321323132 & 12321312132 & & 12313231213 & 23132123213 & 12321231213\\
12313231232 & 12312131232 & 13121323132 & & 12313231232 & 13121323132 & 12312131232\\
12321231213 & 12313231213 & 23132123213 & & 12321231213 & 23132123213 & 12313231213\\
12321231323 & 12132123213 & 12132313213 & & 12321231323 & 12132313213 & 12132123213\\
12321312132 & 12313231213 & 12321323132 & & 12321312132 & 12321323132 & 12313231213\\
12321312132 & 31213123132 & 31232123132 & & 12321312132 & 31232123132 & 31213123132\\
12321323132 & 12313231213 & 12321312132 & & 12321323132 & 12321312132 & 12313231213\\
13121323132 & 12312131232 & 12313231232 & & 13121323132 & 12313231232 & 12312131232\\
13123132312 & 13123212312 & 13231321232 & & 13123132312 & 13231321232 & 13123212312\\
13123212312 & 13123132312 & 13231321232 & & 13123212312 & 13231321232 & 13123132312\\
13212321312 & 13231213123 & 13231232123 & & 13212321312 & 13231232123 & 13231213123\\
13212321312 & 13231321312 & 32123132312 & & 13212321312 & 32123132312 & 13231321312\\
13212321323 & 12131232123 & 13213121323 & & 13212321323 & 13213121323 & 12131232123\\
13213121323 & 12131232123 & 13212321323 & & 13213121323 & 13212321323 & 12131232123\\
13231213123 & 13212321312 & 13231232123 & & 13231213123 & 13231232123 & 13212321312\\
13231213123 & 21312132123 & 21323132123 & & 13231213123 & 21323132123 & 21312132123\\
13231232123 & 13212321312 & 13231213123 & & 13231232123 & 13231213123 & 13212321312\\
13231321232 & 13123132312 & 13123212312 & & 13231321232 & 13123212312 & 13123132312\\
13231321312 & 13212321312 & 32123132312 & & 13231321312 & 32123132312 & 13212321312\\
21231213123 & 21231323123 & 21312132313 & & 21231213123 & 21312132313 & 21231323123\\
21231323123 & 21231213123 & 21312132313 & & 21231323123 & 21312132313 & 21231213123\\
21232131213 & 23121312313 & 23123212313 & & 21232131213 & 23123212313 & 23121312313\\
21312132123 & 13231213123 & 21323132123 & & 21312132123 & 21323132123 & 13231213123\\
21312132313 & 21231213123 & 21231323123 & & 21312132313 & 21231323123 & 21231213123\\
21312313231 & 21312321231 & 21323132123 & & 21312313231 & 21323132123 & 21312321231\\
21312321231 & 21312313231 & 21323132123 & & 21312321231 & 21323132123 & 21312313231\\
21312321231 & 32123213231 & 32131213231 & & 21312321231 & 32131213231 & 32123213231\\
21321232131 & 21323132131 & 23212313231 & & 21321232131 & 23212313231 & 21323132131\\
21323132123 & 13231213123 & 21312132123 & & 21323132123 & 21312132123 & 13231213123\\
21323132123 & 21312313231 & 21312321231 & & 21323132123 & 21312321231 & 21312313231\\
21323132131 & 21321232131 & 23212313231 & & 21323132131 & 23212313231 & 21321232131\\
23121312313 & 21232131213 & 23123212313 & & 23121312313 & 23123212313 & 21232131213\\
23121312321 & 23132123213 & 23132131213 & & 23121312321 & 23132131213 & 23132123213\\
23121312321 & 23132312321 & 31213231321 & & 23121312321 & 31213231321 & 23132312321\\
23123212313 & 21232131213 & 23121312313 & & 23123212313 & 23121312313 & 21232131213\\
23132123213 & 12313231213 & 12321231213 & & 23132123213 & 12321231213 & 12313231213\\
23132123213 & 23121312321 & 23132131213 & & 23132123213 & 23132131213 & 23121312321\\
23132131213 & 23121312321 & 23132123213 & & 23132131213 & 23132123213 & 23121312321\\
23132312131 & 23213121321 & 23213231321 & & 23132312131 & 23213231321 & 23213121321\\
23132312321 & 23121312321 & 31213231321 & & 23132312321 & 31213231321 & 23121312321\\
23212313231 & 21321232131 & 21323132131 & & 23212313231 & 21323132131 & 21321232131\\
23213121321 & 23132312131 & 23213231321 & & 23213121321 & 23213231321 & 23132312131\\
23213231321 & 23132312131 & 23213121321 & & 23213231321 & 23213121321 & 23132312131\\
\hline
\end{tabular}}

\vfill\eject

\centerline{\small\tt
\begin{tabular}{|c|c|c|c|c|c|c|}
\hline
$\f(1)$ & $\f(2)$ & $\f(3)$ & & $\f(1)$ & $\f(2)$ & $\f(3)$ \\
\hline
31213123132 & 12321312132 & 31232123132 & & 31213123132 & 31232123132 & 12321312132\\
31213123212 & 31321232132 & 31321312132 & & 31213123212 & 31321312132 & 31321232132\\
31213212321 & 31213231321 & 31232123132 & & 31213212321 & 31232123132 & 31213231321\\
31213231321 & 23121312321 & 23132312321 & & 31213231321 & 23132312321 & 23121312321\\
31213231321 & 31213212321 & 31232123132 & & 31213231321 & 31232123132 & 31213212321\\
31231323121 & 31232123121 & 32313212321 & & 31231323121 & 32313212321 & 31232123121\\
31232123121 & 31231323121 & 32313212321 & & 31232123121 & 32313212321 & 31231323121\\
31232123132 & 12321312132 & 31213123132 & & 31232123132 & 31213123132 & 12321312132\\
31232123132 & 31213212321 & 31213231321 & & 31232123132 & 31213231321 & 31213212321\\
31321232132 & 31213123212 & 31321312132 & & 31321232132 & 31321312132 & 31213123212\\
31321312132 & 31213123212 & 31321232132 & & 31321312132 & 31321232132 & 31213123212\\
31323121312 & 32131213212 & 32132313212 & & 31323121312 & 32132313212 & 32131213212\\
32123121312 & 32123132312 & 32131213231 & & 32123121312 & 32131213231 & 32123132312\\
32123132312 & 13212321312 & 13231321312 & & 32123132312 & 13231321312 & 13212321312\\
32123132312 & 32123121312 & 32131213231 & & 32123132312 & 32131213231 & 32123121312\\
32123213121 & 32312131231 & 32312321231 & & 32123213121 & 32312321231 & 32312131231\\
32123213231 & 21312321231 & 32131213231 & & 32123213231 & 32131213231 & 21312321231\\
32131213212 & 31323121312 & 32132313212 & & 32131213212 & 32132313212 & 31323121312\\
32131213231 & 21312321231 & 32123213231 & & 32131213231 & 32123121312 & 32123132312\\
32131213231 & 32123132312 & 32123121312 & & 32131213231 & 32123213231 & 21312321231\\
32132313212 & 31323121312 & 32131213212 & & 32132313212 & 32131213212 & 31323121312\\
32312131231 & 32123213121 & 32312321231 & & 32312131231 & 32312321231 & 32123213121\\
32312321231 & 32123213121 & 32312131231 & & 32312321231 & 32312131231 & 32123213121\\
32313212321 & 31231323121 & 31232123121 & & 32313212321 & 31232123121 & 31231323121\\
\hline
\end{tabular}}

\end{document}